\documentclass[twoside]{amsart}
\usepackage{amssymb}
\usepackage{amsmath}
\usepackage{amsthm}
\usepackage[all]{xy}

\newtheorem{theorem}{Theorem}

\theoremstyle{definition}

\theoremstyle{remark}

\newcommand{\cO}{{\mathcal O}}

\newcommand{\cL}{{\mathcal L}}

\newcommand{\cC}{{\mathcal C}}
\newcommand{\cP}{{\mathcal P}}
\newcommand{\cD}{{\mathcal D}}

\newcommand{\C}{\mathbb C}
\newcommand{\R}{\mathbb R}

\newcommand{\Q}{\mathbb Q}

\newcommand{\ra}{\rightarrow}

\begin{document}
\title{Rigidity of Mori cone for Fano manifolds}
\author{Jaros{\l}aw A. Wi\'sniewski}

\thanks{Research supported by a grant of Polish
MNiSzW (N N201 2653 33).}
\address{Instytut Matematyki UW, Banacha 2, PL-02097
Warszawa} \email{J.Wisniewski@mimuw.edu.pl}

\begin{abstract}
Mori cone is rigid in smooth connected families of Fano manifolds.
\end{abstract}

\maketitle

The aim of this note is to give a positive answer to a question raised
at a workshop {\sl Rational curves on Algebraic Varieties} organized
in American Institute of Mathematics, Palo Alto, CA, in May 2007,
\cite[Question 0.7]{PaloAlto}. 

We consider algebraic varieties defined over the field of complex
numbers $\C$. Manifolds are smooth connected varieties. 

If $X$ is a complex projective manifold then by $N^1(X)\subset
H^2(X,\R)$ and $N_1(X)\subset H_2(X,\R)$ we denote $\R$-linear
subspaces spanned by cohomology and homology classes of, respectively,
Cartier divisors and algebraic curves on $X$.  The cone of curves, or
Mori cone of $X$, ${\mathcal C}(X)\subset N_1(X)$ and the cone of nef
divisors ${\mathcal P}(X)\subset N^1(X)$ are $\R^*_{>0}$-spanned by,
respectively, the classes of curves, or effective 1-cycles, and
numerically effective divisors, hence ${\mathcal P}(X):= \{\chi\in
N^1(X): \forall \alpha\in{\mathcal C}(X) \ \ \chi\cdot\alpha\geq 0\}$.
That is, ${\mathcal P}(X)={\mathcal C}(X)^\vee$ in the sense of the
intersection pairing of $N^1(X)$ and $N_1(X)$. See \cite{Mori2} for
more information on these objects.

A manifold $X$ is Fano if its anticanonical divisor $-K_X$ is
ample. If $X$ is Fano then $H^1(\cO_X)=0$ hence $N^1(X)=H^2(X,\R)$ and
$N_1(X)=H_2(X,\R)$.

A smooth family of projective manifolds $\pi: X\ra S$ is, by
definition, a smooth projective morphism of connected complex
manifolds with connected fibers. Geometric fibers of $\pi$ are denoted
by $X_s=\pi^{-1}(s)$ for $s\in S$. By a known result in differential
geometry any such family is topologically locally trivial in complex
topology. A local topological trivialization of $\pi$ induces
identification of homology and cohomology of neighboring fibers.

\begin{theorem}\label{MoriCone}
Let us assume that $\pi:X\ra S$ is a smooth family of Fano manifolds,
that is the relative anticanonical divisor $-K_{X/S}$ is $\pi$-ample.
Then for any two $s_0,\ s_1\in S$ the local identification
$H_2(X_{s_0},\R)=H_2(X_{s_1},\R)$ and
$H^2(X_{s_0},\R)=H^2(X_{s_1},\R)$ yields
$\cC(X_{s_0})=\cC(X_{s_1})$ and $\cP(X_{s_0})=\cP(X_{s_1})$

\end{theorem}

The above result is an immediate consequence of the following main
theorem from \cite[Theorem 1.7]{Duke}, see also \cite{Asian} where the
ideas of \cite{Duke}, based on playing Hard Lefschetz Theorem against
Mori's theory of rational curves, see \cite{Mori1}, \cite{Mori2}, were
further developed. The notation in the present paper is consistent
with that of \cite{Duke}. Recall that given an ample line bundle $\cL$
on a manifold $X$ its nef value, or nef threshold, $\tau(\cL)$ is the
infimum of $t\in\Q$ such that $K_X+t\cL$ is ample $\Q$-divisor. Note
that $\tau(\cL)$ is positive if and only if $K_X$ is not nef.

\begin{theorem}\label{NefValue}
Let us assume that $\pi:X\ra S$ is a smooth (connected) family of
projective manifolds with $\cL$ a $\pi$-ample line bundle. If
$\tau(\cL_{s_0})$ is positive for some $s_0\in S$ then the function
$S\ni s \mapsto \tau(\cL_s)$ is constant.
\end{theorem}

In view of Theorem \ref{NefValue} a proof of Theorem \ref{MoriCone}
can be reduced to choosing appropriate $\pi$-ample line bundle $\cL$.
Namely, suppose that $D_0$ is an ample divisor on $X_{s_0}$. Now,
possibly shrinking $S$ to a connected base over which the family $\pi$
is topologically trivial we can extend $D_0$ to a divisor $\cD$ over
$X$; this extension provides us with an identification
$H^2(X_{s_0},\R)=H^2(X_{s_1},\R)$.  We claim that for a sufficiently
small positive $\epsilon$ the divisor $\cL=-K_X+\epsilon\cdot(\cD +
K_X)$ is $\pi$-ample. Indeed, since ampleness is an open condition,
for any $s\in S$ there exists $\epsilon_s> 0$ such that
$-K_X+\epsilon_s\cdot(\cD + K_X)$ is $\pi$-ample in an open Zariski
neighborhood of $X_s$. Thus we can take $\epsilon={\rm
  min}(\epsilon_{s_0},\epsilon_{s_1})$ and possibly shrink $S$ to a
smaller connected variety containing both $s_0$ and $s_1$. We denote
the resulting restricted family as before: $\pi: X\ra S$. Now we apply
Theorem \ref{NefValue} to argue that $\cD_{s_1}$ is ample: see the
picture below which presents a plane in $N^1(X_s)$ spanned on $K_X$
and $\cD$; it is apparent that invariance of $\tau_s$ yields ampleness
of $\cD_s$ because the line determined by $K_X+\tau_sL_s$ is fixed and
it constitues a border of the nef cone for all $s$. This implies
that, in terms of the identification
$H^2(X_{s_0},\R)=H^2(X_{s_1},\R)$, we get $\cP(X_{s_0})=\cP(X_{s_1})$.
The equality $\cC(X_{s_0})=\cC(X_{s_1})$ follows by duality. This
concludes a proof of Theorem \ref{MoriCone}.
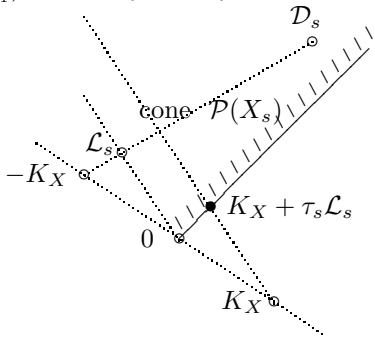
$$
\begin{xy}<12pt,0pt>:
(-4.5,3) ; (4.5,-3) **@{.},
(0,0)*={\circ}="0" , (-1,0)*={0}, 
(-3,2)*={\circ}="1", (-4.5,2)*={-K_X},
(3,-2)*={\circ}="2", (2,-2)*={K_X},
(4.2,6.2)*={\circ}="3", (4,7)*={\cD_s}, 
"1" ; "3" **@{.},
(-1.8,2.7)*={\circ}, (-2.5,3)*={\cL_s},
"0" ; (-3,4.5) **@{.},
"2" ; (-3,7) **@{.},
(1,1)*={\bullet}, (3.5,1)*={K_X+\tau_s\cL_s},
"0" ; (6,6) **@{-},
(1,4)*={{\rm cone}\ \ \cP(X_s)},
(0,0.5) ; (6,6.5) **@{/},
\end{xy}
$$ 

I would like to thank Rob Lazarsfeld who convinced me that it is
worthwhile to write this note. He also observed that rigidity of the
cone of big divisors for Fano manifolds follows from an extension
theorem by Siu, \cite{Siu}, and suggested the following corollary to
Theorem \ref{NefValue}.

\begin{theorem}\label{MinManifolds}
Let $\pi: X\ra S$ be a smooth family of projective manifolds. If for
some $s_0\in S$ the variety $X_{s_0}$ is not minimal, that is
$K_{s_0}$ is not nef, then no variety $X_s$ is minimal.
\end{theorem}

A version of the above theorem was proved in analytic category in
\cite{AndreattaPeternell}.


\end{document}